\pdfoutput=1
\RequirePackage{ifpdf}
\ifpdf 
\documentclass[pdftex]{sigma}
\else
\documentclass{sigma}
\fi

\usepackage{mathdots}
\usepackage{tikz}
\usetikzlibrary{arrows.meta,arrows,calc,decorations.pathmorphing,decorations.pathreplacing,positioning,shapes}
\usepackage{enumitem}

\numberwithin{equation}{section}
\newtheorem{Theorem}{Theorem}[section]

\theoremstyle{definition}
\newtheorem{Definition}[Theorem]{Definition}
\newtheorem{Remark}[Theorem]{Remark}

\newcommand{\Z}{\mathbb{Z}}

\newcommand{\C}{\mathbb{C}}

\begin{document}

\allowdisplaybreaks

\newcommand{\arXivNumber}{2410.21879}

\renewcommand{\thefootnote}{}

\renewcommand{\PaperNumber}{042}

\FirstPageHeading

\ShortArticleName{A New Grounded Partition Identity of Type $D_4^{(3)}$}

\ArticleName{A New Grounded Partition Identity of Type $\boldsymbol{D_4^{(3)}}$\footnote{This paper is a~contribution to the Special Issue on Recent Advances in Vertex Operator Algebras in honor of James Lepowsky. The~full collection is available at \href{https://sigma-journal.com/Lepowsky.html}{https://sigma-journal.com/Lepowsky.html}}}

\Author{Benedek DOMBOS}

\AuthorNameForHeading{B.~Dombos}

\Address{Section de math\'ematiques, Universit\'e de Gen\`eve, Switzerland}
\Email{\mail{benedek.dombos@unige.ch}}

\ArticleDates{Received December 18, 2025, in final form April 19, 2026; Published online April 30, 2026}

\Abstract{In this paper, we prove a new Rogers--Ramanujan-type identity, involving grounded partitions, by computing a character of the affine Kac--Moody algebra \smash{$D_4^{(3)}$} in two different ways. The product side is derived using Lepowsky's product formula, while the sum side is obtained using perfect crystals with a technique of Dousse and Konan.}


\Classification{05A17; 17B67; 05E10}

\renewcommand{\thefootnote}{\arabic{footnote}}
\setcounter{footnote}{0}

\section{Introduction}

A \emph{partition} of a nonnegative integer $ n $ is a non-increasing sequence $\lambda = (\lambda_1, \lambda_2, \dots, \lambda_k)$ of positive integers such that $n = \lambda_1 + \lambda_2 + \dots + \lambda_k$. The terms $ \lambda_i $ are called the \emph{parts} of the partition~$\lambda$. The integer $n$ is called the \emph{size} of $\lambda$, denoted $|\lambda|$. The empty sequence corresponds to the unique partition of size $0$, called the \emph{empty partition}.

In the theory of integer partitions, it is convenient to abbreviate some finite and infinite products by \emph{$q$-Pochhammer symbols}. For formal variables $a$ and $q$, we define
\begin{align*}
 & (a;q)_k := (1 - a)(1 - aq)\bigl(1 - aq^2\bigr) \cdots \bigl(1 - aq^{k-1}\bigr) \qquad\text{for $ k \geq 1 $},\\
 & (a;q)_0:= 1, \qquad (a;q)_\infty := \prod_{k=0}^{\infty} \bigl(1 - aq^k\bigr).
\end{align*}
We also use the shorthand notation
\[(a_1, \ldots, a_\ell; q^n)_\infty := (a_1; q^n)_\infty \cdots (a_\ell; q^n)_\infty \qquad\text{for $ \ell \geq 1 $}.\]
Partition identities of the form: the number of partitions of $n$ with certain congruence conditions on the parts equals the number of partitions of $n$ with certain difference conditions on the parts, are called identities of \emph{Rogers--Ramanujan type}, named after the following theorem.

\begin{Theorem}[Rogers--Ramanujan identities, in terms of $q$-series, \cite{RR19}]
\begin{align*}
 \sum_{n=0}^{\infty} \frac{q^{n^2}}{(q;q)_n} = \frac{1}{\bigl(q;q^5\bigr)_\infty \bigl(q^4;q^5\bigr)_\infty}, \qquad
 \sum_{n=0}^{\infty} \frac{q^{n^2+n}}{(q;q)_n} = \frac{1}{\bigl(q^2;q^5\bigr)_\infty \bigl(q^3;q^5\bigr)_\infty}.
\end{align*}
\end{Theorem}
Each of these $q$-series identities corresponds to a partition identity. Here, we only state the first identity (see \cite{Mac16} and \cite{Sch17} for further details). The number of partitions of $ n $ such that the difference between consecutive parts is at least 2 is equal to the number of partitions of $ n $ into parts congruent to $ 1 $ or $ 4 \pmod 5 $.

Lepowsky and Wilson~\cite{LW84,LW85} gave a representation-theoretic
proof of the Rogers--Ramanujan identities, building on earlier work of
Lepowsky and Milne~\cite{LM78}. The product side can be obtained
directly from the principal specialisation of the Weyl--Kac character
formula (see Kac~\cite[Theorem 10.4]{Kac90}), for highest
weight characters of type \smash{$A_1^{(1)}$} at level $3$. The
representation-theoretic interpretation of the sum side is more
intricate. In these papers, Lepowsky and Wilson developed and applied
the theory of vertex operators to derive it.

Grounded partitions are combinatorial objects introduced by Dousse and
Konan in~\cite{DK19b} and~\cite{DK22}, motivated by the theory of
perfect crystals, introduced by Kang, Kashiwara, Misra, Miwa,
Nakashima, and Nakayashiki~\cite{KMN92}. Here, we
define only a simplified version of grounded partitions with exact
difference conditions, which will be sufficient for our purposes.

\begin{Definition}\label{def:grounded_partitions}
 Let $ n $ be a nonnegative integer, and let
 $C=\{c_0,c_1,\ldots,c_n\}$ denote the set of \emph{colours}, let $M$ be a square matrix of size $n+1$ with integer entries.
 Fix a colour $c_i\in C$. A~\emph{grounded partition}
 $ \pi = (\pi_1, \pi_2, \dots ) $ corresponding to the
 matrix $M$ with ground $0_{c_i}$ is defined as follows:
\begin{enumerate}[label=(\roman*)]\itemsep=0pt
\item Each part $ \pi_j $ is a positive integer indexed by some colour in $C$.
\item We start with the part $0_{c_i}$, as an initial condition,
 called the \emph{ground}.
\item Exact difference conditions: for $ j \geq 0 $, let
 $\pi_{j}\leq\pi_{j+1}$ be two consecutive parts with
 colours~$c_{i_j}$ and $c_{i_{j+1}}$, respectively. Then
 \[
 \pi_{j+1}-\pi_{j} = (i_{j+1},i_j)\text{-entry in } M.
 \]
\end{enumerate}
The \emph{size} of $\pi$, denoted $|\pi|$, is defined as the
sum of the positive integers obtained by forgetting the colours of the parts. We denote the set of grounded partitions corresponding to the matrix~$M$ with ground $0_{c_i}$ by \smash{$\mathcal{P}_{M}^{c_i}$}.
\end{Definition}

\begin{Remark}
Partitions are typically written in weakly decreasing order. For
grounded partitions, however, we list parts in weakly increasing
order, as the colour of the $ 0 $-part is significant.
\end{Remark}

Dousse and Konan~\cite{DK19b,DK22, DK22c} provided combinatorial character formulas for all level-one representations in affine types $A$, $B$, $C$, $D$ given by generating functions of grounded partitions. Furthermore, Dousse, Hardiman, and Konan~\cite{DHK25} found identities in affine type \smash{$A_1^{(1)}$} for all levels~${\ell \geq 1}$.

\begin{Theorem}[{\cite[Theorem~1.6]{DHK25}}]
Let $M_\ell$ denote the square matrix of size $\ell+1$ whose
$(i,j)$-entry is $|i-j|$. Then, the generating function for grounded partitions in
\smash{$ \mathcal{P}_{M_\ell}^{c_i} $} is given by
\[
 \sum_{\lambda \in \mathcal{P}_{M_\ell}^{c_i}} q^{|\pi|} = \frac{\bigl(q^{i+1}, q^{\ell-i+1}, q^{\ell+2}; q^{\ell+2}\bigr)_\infty}{\bigl(q;q^2\bigr)_\infty (q;q)_\infty}.
\]
\end{Theorem}
They used the character formula for perfect crystals~\cite{KMN92} corresponding to
the standard modules of the affine Lie algebra \smash{$ A_1^{(1)} $} at
level $ \ell $, of highest weight
$ i \Lambda_0 + (\ell-i) \Lambda_1 $. Using the associated energy
matrices, they showed that this character formula coincides with the
generating function for grounded partitions $ \mathcal{P}_{M_\ell, i} $,
and obtained the right-hand side via the principal specialization of
the Weyl--Kac character formula.

In this paper, we prove a new identity of Rogers--Ramanujan-type using a perfect crystal of affine type \smash{$D_4^{(3)}$}.

\begin{Theorem}\label{thm:G2identity}
Let $M_{G_2}$ be the following square matrix with rows and columns indexed by colours~$\{a,b,c,d,e,f,g,h\}$:
\begin{align*}
M_{G_2} = \begin{array}{c|cccccccc}
& a & b & c & d & e & f & g & h \\
\hline
a &8&5&6&3&4&5&6&7\\
b &7&8&5&6&7&4&5&6\\
c &6&7&8&5&6&7&4&5\\
d &5&6&7&4&5&6&3&4\\
e &4&5&6&7&8&5&6&3\\
f &3&4&5&6&7&8&5&2\\
g &2&3&4&5&6&7&8&1\\
h &1&2&3&4&5&6&7&0\\
\end{array}.
\end{align*}
For a non-negative integer $n$, the number of grounded partitions in \smash{$\mathcal{P}_{M_{G_2}}^h$} of size $n$ is equal to the number of partitions of $n$ with parts congruent to $1$ or $5$ modulo $6$.
\end{Theorem}

Dousse and Konan~\cite{DK19b} also considered a variation of Definition~\ref{def:grounded_partitions}. \emph{Grounded partitions with non-exact difference conditions} are defined in the same way as grounded partitions except condition~(iii):
\begin{enumerate}\itemsep=0pt
\item[(iii$'$)] Non-exact difference conditions: for $ j \geq 0 $, let
 $\pi_{j}\leq\pi_{j+1}$ be two consecutive parts with
 colours $c_{i_j}$ and $c_{i_{j+1}}$, respectively. Then
 \[
 \pi_{j+1}-\pi_{j} \geq (i_{j+1},i_j)\text{-entry in } M.
 \]
\end{enumerate}
Let \smash{$\mathcal{P}_{M}^{c_i,\geq}$} denote the set of grounded partitions with non-exact difference conditions corresponding to the matrix $M$ with ground $0_{c_i}$. Furthermore, we denote by~\smash{$\mathcal{P}_{M,m}^{c_i,\geq}$} the set of elements of~\smash{$\mathcal{P}_{M}^{c_i,\geq}$} satisfying the following extra condition:
\begin{enumerate}\itemsep=0pt
\item[(iv)] Congruence conditions: for any $j\in\{1,\ldots,n\}$, parts with colour $c_j$ are congruent to the $ (i,j)$-entry in $M$ modulo $m$.
\end{enumerate}

\begin{Theorem}\label{thm:G2identity_geq}
For a non-negative integer $n$, the number of elements in \smash{$\mathcal{P}_{M_{G_2},4}^{h,\geq}$} of size $n$ is equal to the number of partitions with parts congruent to $1$ or $5$ modulo $6$, or congruent to $0$ modulo $4$.
\end{Theorem}

In Section~\ref{section:preliminaries}, we give a brief summary of key definitions and results on affine Kac--Moody algebras and perfect crystals. In Section~\ref{section:cong}, we derive an explicit infinite-product form for the principal specialisation of the character of the level~1 irreducible representation with highest weight $\Lambda_0$ of the twisted affine Kac--Moody algebra of type \smash{$D_4^{(3)}$}, expressing it as a generating function for partitions with congruence conditions. In Section~\ref{section:difference}, we express this character as a generating function for grounded partitions, using a perfect crystal of type \smash{$D_4^{(3)}$} and the combinatorial framework of Dousse--Konan, which completes the proof of Theorems~\ref{thm:G2identity}~and~\ref{thm:G2identity_geq}. In Section~\ref{section:recursions}, we derive recursions for this generating function and show that it is impossible to keep track of any colours in the generating functions and still have an infinite product expression.

\section{Preliminaries} \label{section:preliminaries}

In this section, we provide a summary of key constructions and results on affine Kac--Moody algebras and perfect crystals. For an introduction to affine Kac--Moody algebras, two standard references are Kac~\cite{Kac90} and Carter~\cite{Car10}. For an introduction to perfect crystals, a standard reference is Hong--Kang~\cite[Chapter~10]{HK02}.

\subsection{Lepowsky's product formula}

The definition of Kac--Moody algebras is motivated by Serre's reconstruction of finite-dimen\-sional semisimple Lie algebras from their Cartan matrices. A \emph{generalised Cartan matrix} is a~square matrix $A = (a_{ij})_{i,j \in I}$, where $I = \{0, 1, \ldots, n\}$, with integer entries such that $a_{ii} = 2$ for all $i \in I$; $a_{ij} \leq 0$ for all $i \neq j$; and $a_{ij} = 0$ if and only if $a_{ji} = 0$. Given a generalised Cartan matrix $A = (a_{ij})_{i,j \in I}$, the associated \emph{Kac--Moody algebra} is defined as the Lie algebra~$\mathfrak{g}(A)$ over~$\mathbb{C}$, generated by the elements $\{e_i, h_i, f_i\}_{i \in I}$, subject to the following relations:
\begin{enumerate}[label=(\roman*)]\itemsep=0pt
 \item $[h_i,h_i] =0$ for all $i \in I$,
 \item $[e_i,f_i] =h_i$ for all $i \in I$,
 \item $[e_i,f_j] =0$ for all $i \neq j$,
 \item $[h_i,e_j] =a_{ij}e_j$, $[h_i,f_j] =-a_{ij}f_j$ for all $i,j \in I$,
 \item \smash{$\mathrm{ad}_{e_i}^{1-a_{ij}}e_j = \mathrm{ad}_{f_i}^{1-a_{ij}}f_j = 0$} if $i \neq j$, where $\mathrm{ad}$ is the \textit{adjoint}.
\end{enumerate}

A generalised Cartan matrix $A$ is called \emph{symmetrisable} if there exist $d_i \in \Z$ for all $i \in I$, such that $(a_{ij}d_i)_{i,j \in I}$ is a symmetric matrix. If a symmetrisable generalised Cartan matrix $A$ has nonzero determinant, it is said to be of \emph{finite type}. If $A$ has determinant zero and each proper principal minor is positive definite, it is said to be of \emph{affine type}. A generalised Cartan matrix of finite type and size $n+1$ has rank $n+1$, whereas a generalised Cartan matrix of affine type and size $n+1$ has rank $n$. For a generalised Cartan matrix $A$ of finite type, $\mathfrak{g}(A)$ is a finite dimensional semisimple Lie algebra, and all finite dimensional semisimple Lie algebras arise this way. For a generalised Cartan matrix $A$ of affine type, $\mathfrak{g}(A)$ is called a \emph{Kac--Moody algebra of affine type}.

For the remainder of this subsection, let $\mathfrak{g}$ be a Kac--Moody algebra of affine type. The \emph{Kac labels} $a_i$, $i\in I$ are defined as the unique coprime positive integers satisfying
\[
\sum_{j\in I} a_{ij} a_j = 0 \qquad\text{for all } i\in I.
\]
The \emph{dual Kac labels} $a_i^\vee$, $i\in I$, are defined as the unique coprime positive integers satisfying
\[
\sum_{i\in I} a_i^\vee a_{ij} = 0 \qquad\text{for all } j\in I.
\]
The \emph{central element} $c\in\mathfrak{g}$ is defined as
\[
c := \sum_{i\in I} a_i^\vee h_i,
\]
and the \emph{degree derivation} $d$ is defined as the unique element of the Kac--Moody algebra $\mathfrak g$ satisfying ${[d,e_0]=e_0}$, ${[d,f_0]=-f_0}$, ${[d,e_i]=[d,f_i]=0}$ for all $i\neq 0$, and $[d,h_i]=0$ for all $i\in I$. The \emph{Cartan subalgebra} $\mathfrak{h}$ of~$\mathfrak{g}$ (that we choose) is the maximal abelian subalgebra given by $\mathfrak{h} := \C \{h_i \mid i\in I\} \oplus \C d$.
A~$ \mathfrak g$-module $V$ is called a \emph{highest weight module} of highest weight $\Lambda\in\mathfrak h^{*}$ if there exists a non-zero vector~$v_\Lambda\in V$ such that $h.v_\Lambda = \Lambda(h)v_\Lambda$ for all $h\in \mathfrak{h}$; $e_i.v_\Lambda = 0$ for all $i \in I$; and $V$ is generated by $v_\Lambda$ under the action of the operators $f_i$, $i\in I$. A highest weight module $V$ admits a unique highest weight $\Lambda$ up to scaling and a \emph{weight-space decomposition}:
\[
V = \bigoplus_{\lambda \in \mathfrak{h}^*} V_\lambda,\qquad V_\lambda := \{v \in V \mid h.v = \lambda(h)v\ \forall h \in \mathfrak{h} \}.
\]
The \emph{character} of $V$ is defined as
\[\operatorname{ch} V := \sum_{\lambda \in \mathfrak{h}^*}(\dim V_\lambda) {\rm e}^\lambda.\]
The set of \emph{dominant integral weights} of $\mathfrak{g}$ is defined as $P_{+} := \{\Lambda\in\mathfrak h^{*} \mid \Lambda(h_i)\in\mathbb{Z}_{\ge0}\ \text{for all } i\in I\}$, and the \emph{weight lattice} is defined as $P := \{\Lambda\in\mathfrak h^{*} \mid \Lambda(h_i)\in\mathbb{Z}\ \text{for all } i\in I\}$. We only consider highest weight modules of $\mathfrak{g}$ with highest weight $\Lambda\in P^+$.

The elements $h_i$, $i\in I$, are called \emph{simple coroots}, and the \emph{simple roots} $\alpha_i\in\mathfrak h^*$ are defined by $\alpha_i(h_j)=a_{ji}$ for $j\in I$, $\alpha_i(c)=0$ and $\alpha_i(d)=\delta_{i0}$. For each $i\in I$, consider the reflection $s_i\in\operatorname{Aut}_{\mathbb C}(\mathfrak h)$ by $s_i(h_j)=h_j-a_{ji}h_i$, $j\in I$, and $s_i(d)=d$. The \emph{Weyl group} is defined as $W:=\langle s_i \mid i\in I \rangle$.
For any $w\in W$, we denote by $|w|$ the length of $w$ in $W$ with respect to the generators $s_i$, $i\in I$. The Weyl group $W$ acts on $\mathfrak h^{*}$ by $(w\cdot\lambda)(h)=\lambda\bigl(w^{-1}h\bigr)$, where $w\in W$, $\lambda\in\mathfrak h^{*}$, and $h\in\mathfrak h$.

The set of \emph{real roots} of $\mathfrak{g}$ is defined as $\Phi_{\mathrm{Re}}:= W\{\alpha_i \mid i\in I\}$, the set of \emph{imaginary roots} is defined as $\Phi_{\mathrm{Im}}:=\{n\delta \mid n\in\mathbb Z\setminus\{0\}\}$, and the set of \emph{roots} of $\mathfrak{g}$ is defined as $\Phi:=\Phi_{\mathrm{Re}}\cup\Phi_{\mathrm{Im}}$. We have the \emph{root space decomposition}
\[
\mathfrak{g} = \mathfrak{h}  \oplus  \bigoplus_{\alpha \in \Phi} \mathfrak{g}_{\alpha},
\qquad
\mathfrak{g}_{\alpha} := \{ x \in \mathfrak{g} \mid [h,x] = \alpha(h).x \ \text{for all } h \in \mathfrak{h} \}.
\]
The set of \emph{positive roots}, denoted $\Phi_+$, consists of roots $\alpha\in\Phi$ that can be written as a nonnegative linear combination of simple roots. The \emph{basic imaginary root} is defined by
\[
\delta:=\sum_{i\in I} a_i \alpha_i.
\]
Let $A_{\mathrm{fin}}$ denote the generalised Cartan matrix obtained by deleting the $0$-th row and column of $A$, and let $\mathfrak{g}_{\mathrm{fin}}$ be the finite-dimensional Lie algebra corresponding to $A_{\mathrm{fin}}$.
Let $\Phi_{\mathrm{fin}}$ denote the root system of $\mathfrak{g}_{\mathrm{fin}}$. Then the real roots of $\mathfrak g$ can be parametrised as $\Phi_{\mathrm{Re}}
 = \{ \alpha + n\delta \mid \alpha\in \Phi_{\mathrm{fin}},\, n\in\mathbb Z \}$.

 The \emph{level} of an element $\lambda \in \mathfrak{h}^*$ is defined as $\mathrm{lev}(\lambda) := \lambda(c)$. The \emph{fundamental weights} $\Lambda_i\in\mathfrak h^{*}$, $i\in I$ are defined by $\Lambda_i(h_j)=\delta_{ij}$ for $j\in I$ and $\Lambda_i(d)=0$. We have
 \begin{align*}
 \Lambda_i(c)=\sum_{j\in I}a_j^\vee\Lambda_i(h_j)=a_i^\vee.
 \end{align*}
 The \emph{Weyl vector} $\rho\in\mathfrak h^{*}$ is defined by $\rho(h_i)=1$ for $i\in I$. Let $\Lambda=\sum_{i\in I} m_i\Lambda_i$ with ${m_i\in\mathbb Z_{\ge0}}$.
Then $\Lambda\in P_{+}$, and there exists a unique (up to isomorphism) irreducible highest-weight
$\mathfrak g$\nobreakdash-module~$L(\Lambda)$ of highest weight $\Lambda$. Its level is given by
\[
\operatorname{lev}(\Lambda)=\sum_{i\in I} m_i\Lambda_i(c)=\sum_{i\in I}m_i\sum_{j\in I}a_j^\vee\Lambda_i(h_j)=\sum_{i\in I}a_i^\vee m_i.
\]
We denote by $P_\ell^+$ the set of dominant integral weights of level~$\ell$.

The Weyl--Kac character formula gives an explicit expression for characters of highest weight modules of $\mathfrak{g}$ with dominant integral highest weight in terms of the action of the Weyl group on~roots.
\begin{Theorem}[{Weyl--Kac character formula~\cite[Theorem~10.4]{Kac90}}]
The character of the irreducible representation $L(\lambda)$ of $\mathfrak{g}$ with highest weight $ \Lambda \in P_+$ can be expressed as
\begin{align*}
\operatorname{ch}L(\Lambda) = \frac{\sum_{w \in W} (-1)^{|w|} {\rm e}^{w(\lambda + \rho) - \rho}}{\prod_{\alpha \in \Phi^+} (1 - {\rm e}^{-\alpha})^{m(\alpha)}},
\end{align*}
where $m(\alpha)$ denotes $\dim (\mathfrak{g}_{\alpha})$.
\end{Theorem}

There is a specialisation of the Weyl--Kac character formula that yields a product expansion. This expansion can be used to obtain the product side of $q$-series identities, corresponding to congruence conditions in partition identities. For a sequence $s = (s_0,s_1,\ldots,s_l)$ of non-negative integers, the $s$-\emph{specialisation} is defined as
\[
\mathrm{F}_s\colon \ \Z[[{\rm e}^{-\alpha_0},{\rm e}^{-\alpha_1},\ldots,{\rm e}^{-\alpha_l}]]\to \Z[[q]],\qquad {\rm e}^{-\alpha_i}\mapsto q^{s_i}.
\]
The special case $s=(1,1,\ldots,1)$ is called the \emph{principal specialisation} and is denoted by $\mathrm{F}_1$. Consider the product
\[
D(\Phi):=\prod_{\alpha\in\Phi^+} (1-{\rm e}^{-\alpha} )^{m(\alpha)}.
\]
Then the principal specialisation of the character $\operatorname{ch}L(\Lambda)$ of highest weight $\Lambda$ has the following product expansion.

\begin{Theorem}[Lepowsky's product formula, \cite{Lep78}]
\begin{align*}
 \mathrm{F}_1\bigl({\rm e}^{-\Lambda}\operatorname{ch}L(\Lambda)\bigr) = \frac{\mathrm{F}_{(\Lambda(h_0)+1,\ldots,\Lambda(h_l)+1)} D\bigl(\Phi^\vee\bigr)}{\mathrm{F}_1 D\bigl(\Phi^\vee\bigr)},
\end{align*}
where $\Phi^\vee$ denotes the dual root system.
\end{Theorem}
\begin{Remark}
In the denominator, it makes no difference whether we use the root system or its dual, as $\mathrm{F}_1 D\bigl(\Phi^\vee\bigr) = \mathrm{F}_1 D(\Phi)$, see \cite[equation~(10.8.5)]{Kac90}.
\end{Remark}

\subsection{Perfect crystals}
Let $\mathfrak{g}$ be a Kac--Moody algebra with weight lattice $P$, simple roots $\{\alpha_i\mid i\in I\}$ and coroots~${\{h_i \mid i\in I\}}$.
\begin{Definition}
A \emph{crystal} $ \mathcal{B} $ is a nonempty set endowed with the maps: \smash{$ {\rm wt}\colon \mathcal{B} \to P $} and \smash{$ \tilde{e}_i,\tilde{f}_i\colon \mathcal{B} \to \mathcal{B} \cup \{ 0 \},$} are called \emph{raising} and \emph{lowering operators}, respectively, satisfying the following conditions for all $ i\in I $:
\begin{enumerate}[label=(\roman*)]\itemsep=0pt
 \item if $ \tilde{f}_i(v) \neq 0 $, then $ {\rm wt}\bigl(\widetilde{f}_i(v)\bigr) = {\rm wt}(v) - \alpha_i $,
 \item if $ \tilde{e}_i(v) \neq 0 $, then $ {\rm wt} (\widetilde{e}_i(v) ) = {\rm wt}(v) + \alpha_i $,
 \item $\tilde{f}_i(v)=w$ if and only if $\tilde{e}_i(w)=v$ for all $v,w\in \mathcal{B}$,
 \item $ \varphi_i(v) = \varepsilon_i(v) + ({\rm wt}(v) ) (h_i )$, where
 \begin{itemize}\itemsep=0pt
 \item $ \varepsilon_i(v) := \max \bigl\{ k \geq 0 \mid \tilde{e}_i^k(v) \neq 0 \bigr\} $,
 \item $ \varphi_i(v) := \max \bigl\{ k \geq 0 \mid \tilde{f}_i^k(v) \neq 0 \bigr\} $.
 \end{itemize}
\end{enumerate}
\end{Definition}

The \emph{crystal graph} associated to $\mathcal{B}$ has vertex set $\mathcal{B}$, and oriented, coloured edges
\[
v_1 \overset{i}{\longrightarrow} v_2 \qquad \text{if and only if} \qquad \tilde{f}_i(v_1) = v_2 \qquad \text{(or equivalently } \tilde{e}_i(v_2) = v_1).
\]
When the operator \smash{$\tilde{f}_i$} takes a vertex to $0$, we simply omit the arrow.

\begin{Definition}\label{def:crystaltensor}
Given two crystals $ \mathcal{B}_1 $ and $ \mathcal{B}_2 $, the \emph{tensor product} of crystals $ \mathcal{B}_1 \otimes \mathcal{B}_2 $ is a new crystal defined on the set $ \mathcal{B}_1 \times \mathcal{B}_2 $. The crystal structure on $ \mathcal{B}_1 \otimes \mathcal{B}_2 $ is given by the following rules for elements $ v_1 \in \mathcal{B}_1 $ and $ v_2 \in \mathcal{B}_2 $:
\[
{\rm wt}(b_1 \otimes b_2) := {\rm wt}(v_1) + {\rm wt}(v_2),
\]
\[
\tilde{f}_i(v_1 \otimes v_2) :=
\begin{cases}
\tilde{f}_i(v_1) \otimes v_2 & \text{if } \varphi_i(v_1) > \varepsilon_i(v_2), \\
v_1 \otimes \tilde{f}_i(v_2) & \text{if } \varphi_i(v_1) \leq \varepsilon_i(v_2).
\end{cases}
\]
\end{Definition}

To each highest weight $U_q(\mathfrak{g})$-module of highest weight $\Lambda$, one can associate a unique crystal of highest weight $\Lambda$. See Hong--Kang~\cite[Chapters~3 and~4]{HK02} for the definition of the quantum affine algebra $U_q(\mathfrak{g})$ and the construction of crystal bases. Perfect crystals are associated with finite-dimensional modules over the derived quantum algebra $U_q'(\mathfrak{g})$. The following definition is due to Kang, Kashiwara, Misra, Miwa, Nakashima, and Nakayashiki~\cite{KMN92}.

\begin{Definition}\label{def:perfectcrystal} For a positive integer $\ell$, a crystal $\mathcal{B}$ with finitely many vertices is said to be a~\emph{perfect crystal} of level $\ell$ for $\mathfrak{g}$ if the following conditions hold:
\begin{enumerate}[label=(\roman*)]\itemsep=0pt
 \item $\mathcal{B} \otimes \mathcal{B}$ is connected.
 \item There exists $\lambda_0\in P$ such that
 \begin{align*}
 {\rm wt}(\mathcal{B}) \subset \lambda_0 + \frac{1}{a_0} \sum_{i \neq 0} \mathbb{Z}_{\leq 0} \alpha_i
 \qquad \text{and} \qquad
 |\mathcal{B}_{\lambda_0}| = 1,
 \end{align*}
 where $a_0$ is the coefficient of the simple root $\alpha_0$ in the basic imaginary root $\delta$ and $\mathcal{B}_{\lambda_0}$ denotes the set of vertices in $\mathcal{B}$ of weight $\lambda_0$.
 \item For any $v \in \mathcal{B}$, we have $ (\varepsilon(v) )(c) \geq \ell$, where $c$ is the canonical central element and the operator $\varepsilon$ on $\mathcal{B}$ is given by
 \begin{align*}
 \varepsilon(v) := \sum_{i=0}^n \varepsilon_i(v) \Lambda_i,
 \end{align*}
 where $\Lambda_i$, $i\in I$, are the basic fundamental weights of $\mathfrak{g}$.
 \item For every $\lambda\in P_\ell^+$, there exist unique vertices $v^\lambda$ and $v_\lambda$ in $\mathcal{B}$ such that $\varepsilon\bigl(v^\lambda\bigr) = \lambda$ and $\varphi(v_\lambda) = \lambda$, where the operator $\varphi$ is defined as
 \[\varphi(v) := \sum_{i=0}^n \varphi_i(v) \Lambda_i.\]
\end{enumerate}
\end{Definition}

\begin{Definition}\label{def:groundstatepath}
For $\lambda \in P^+_\ell$, the \emph{ground state path of weight} $\lambda$ is the tensor product
\[
\mathbf{p}_\lambda = (\phi_k)_{k=0}^{\infty} = \cdots \otimes \phi_{k+1} \otimes \phi_k \otimes \cdots \otimes \phi_1 \otimes \phi_0,
\]
where $\phi_0$ is the unique vertex in $b_{\lambda_0}\in\mathcal{B}$ of weight $\lambda_0:= \lambda$ and
\[
\lambda_{k+1} := \varepsilon\bigl(b_{\lambda_k}\bigr), \qquad \phi_{k+1} := b_{\lambda_{k+1}} \qquad \text{for all } k \geq 0.
\]
A tensor product $\mathbf{p} = (p_k)_{k=0}^{\infty} = \cdots \otimes p_{k+1} \otimes p_k \otimes \cdots \otimes p_1 \otimes p_0$ of elements $p_k \in \mathcal{B}$ is said to be a \emph{$\lambda$-path} if $p_k = \phi_k$ for all $k$ large enough. Let $\mathcal{P}(\lambda)$ denote the set of $\lambda$-paths.
\end{Definition}
The tensor product rule in Definition \ref{def:crystaltensor} endows the set $\mathcal{P}(\lambda)$ with a $\mathfrak{g}$-crystal structure. There is a crystal isomorphism \smash{$\mathcal{B}(\lambda) \overset{\sim}{\longrightarrow} \mathcal{P}(\lambda)$}, see \cite[Theorem 10.6.4]{HK02}, where $\mathcal{B}(\lambda)$ denotes the affine crystal corresponding to the irreducible $\mathfrak{g}$-module of highest weight $\lambda$.
\begin{Definition}\label{def:energy}
An \emph{energy function} on $\mathcal{B} \otimes \mathcal{B}$ is a map $H\colon \mathcal{B} \otimes \mathcal{B} \rightarrow \mathbb{Z}$ such that for all $i \in \{0, \dots, n-1\}$ and $b_1$, $b_2$ with $f_i(b_1 \otimes b_2) \neq 0$, we have
\begin{align*}
H\left(f_i(b_1 \otimes b_2)\right) &=
\begin{cases}
H(b_1 \otimes b_2) & \text{if } i \neq 0, \\
H(b_1 \otimes b_2) - 1 & \text{if } i = 0 \text{ and } \varphi_0(b_1) > \varepsilon_0(b_2), \\
H(b_1 \otimes b_2) + 1 & \text{if } i = 0 \text{ and } \varphi_0(b_1) \leq \varepsilon_0(b_2).
\end{cases}
\end{align*}
Given a total order on $\mathcal{B}$, the square matrix $H$ with rows and columns indexed by $\mathcal{B}$ whose \mbox{$(v_1,v_2)$-th} entry is $H(v_1 \otimes v_2)$ is called the \emph{energy matrix} of $\mathcal{B}$.
\end{Definition}

For any perfect crystal, there exists an energy function, which is uniquely determined by a~single value, as $\mathcal{B}\otimes\mathcal{B}$ is connected. The following character formula was also proved in~\cite{KMN92}, and is known as the (KMN)$^2$-character formula.

\begin{Theorem}[(KMN)$^2$-character formula, \cite{KMN92}]\label{thm:KMN2}
Let $\mathcal{B}$ be a perfect crystal for an affine Kac--Moody algebra $\mathfrak{g}$ and $\lambda \in P_\ell^+$ with $\lambda(d)=0$. Let $H$ the an energy function on $\mathcal{B} \otimes \mathcal{B}$ and let $\mathbf{p} = (p_k)_{k=0}^{\infty} \in \mathcal{P}(\lambda)$. Then the character of the irreducible highest weight module $L(\lambda)$ is given by the following expression:
\begin{align*}
\operatorname{ch}L(\lambda)
&= \sum_{p \in \mathcal{P}(\lambda)} {\rm e}^{\operatorname{wt}p},
\end{align*}
where the weight of a path is given by
\begin{align*}
\operatorname{wt}p
&= \lambda + \sum_{k=0}^{\infty} \bigl( \overline{\mathrm{wt}}\,p_k - \overline{\mathrm{wt}}\,\phi_k \bigr) - \frac{\delta}{a_0} \sum_{k=0}^{\infty} (k+1) ( H(p_{k+1} \otimes p_k) - H(\phi_{k+1} \otimes \phi_k) ),
\end{align*}
and $\overline{{\rm wt}}(v)$ is defined as $\varphi(v) - \varepsilon(v)$ modulo $\delta$ such that the coefficient of $\alpha_0$ in $\overline{{\rm wt}}(v)$ is $0$.
\end{Theorem}

Note that eventually we have $p_k = \phi_k$, making this sum finite. A specialisation of Theorem~\ref{thm:KMN2} yields the following formula in the case where the ground state path is constant. Assume $\lambda \in P^+_\ell$ is such that $b_\lambda = b^\lambda = \phi$, and set $H(\phi \otimes \phi) = 0$. Then $\overline{{\rm wt}}\, \phi = 0$, and we have
\begin{align*}
\operatorname{wt}p &= \lambda + \sum_{k=0}^{\infty} \Biggl( \overline{{\rm wt}}\,p_k - \frac{\delta}{d_0} \sum_{j=k}^{\infty} H(p_{j+1} \otimes p_j) \Biggr).
\end{align*}

\section{Congruence conditions}\label{section:cong}

In this section, we derive an explicit infinite product form for the principal specialisation of the character of the irreducible representation of the twisted affine Kac--Moody algebra of type~\smash{$D_4^{(3)}$} of highest weight $\Lambda_0$, expressing it as a generating function for partitions with congruence conditions.

The affine Lie algebras \smash{$G_2^{(1)}$} and \smash{$D_4^{(3)}$} are defined by the generalised Cartan matrices shown in Figure~\ref{fig:CMs}.

\begin{figure}[!t]
 \centering
$\begin{array}{c|ccc}
& {\color{black}0} & {\color{black}1} & {\color{black}2} \\
\hline
 {\color{black}0} & \,\,2 & -1 & \,\,0 \\
 {\color{black}1} & -1 & \,\,2 & -1 \\
 {\color{black}2} & \,\,0 & -3 & \,\,2
\end{array}$ \qquad
$\begin{array}{c|ccc}
& {\color{black}0} & {\color{black}1} & {\color{black}2} \\
\hline
 {\color{black}0} & \,\,2 & -1 & \,\,0 \\
 {\color{black}1} & -1 & \,\,2 & -3 \\
 {\color{black}2} & \,\,0 & -1 & \,\,2
\end{array}$

\vspace{1mm}

(a) $G_2^{(1)}$ \hspace{25mm} (b) $D_4^{(3)}$

\vspace{2mm}

 \caption{Generalised Cartan matrices of \smash{$G_2^{(1)}$} and \smash{$D_4^{(3)}$}.}
 \label{fig:CMs}
\end{figure}

The finite part of both \smash{$G_2^{(1)}$} and \smash{$D_4^{(3)}$} is the simple Lie algebra of type $G_2$, whose root system, with simple roots $\alpha_1$ and $\alpha_2$, is shown in Figure~\ref{fig:G2roots}.
\begin{figure}[!t]
\centering
\begin{tikzpicture}

 \foreach \ang in {60,120,180,240,300,360} {
 \draw[->,black!80!black,thick] (0,0) -- (\ang:2cm);
 }

 \foreach \ang in {30,90,150,210,270,330} {
 \draw[->,black!80!black,thick] (0,0) -- (\ang:3cm);
 }


 \node[scale=0.6, anchor=south east] at (150:3cm) {$\alpha_1$};

 \node[scale=0.6, anchor=south] at (120:2.85cm) {$\alpha_1+\alpha_2$};

 \node[scale=0.6, anchor=south] at (90:3cm) {$2\alpha_1+3\alpha_2$};

 \node[scale=0.6, anchor=west] at (60:2cm) {$\alpha_1+2\alpha_2$};

 \node[scale=0.6, anchor=west] at (30:3cm) {$\alpha_1+3\alpha_2$};

 \node[scale=0.6, anchor=east, xshift=-0.1cm] at (0:3cm) {$\alpha_2$};

\end{tikzpicture}
\caption{$G_2$ root system.}\label{fig:G2roots}
\end{figure}
In particular, the set of short roots of the finite part $G_2$ is given by
\[
\Phi_s = \pm \{\alpha_2, \alpha_1 + \alpha_2, \alpha_1 + 2\alpha_2\}.
\]
and the set of long roots of $G_2$ is given by{\samepage
\[
\Phi_l = \pm \{\alpha_1, 2\alpha_1 + 3\alpha_2, \alpha_1 + 3\alpha_2\},
\]
and the set of roots of $G_2$ is given by $\Phi=\Phi_s\sqcup\Phi_l$.}

Next, we need to parametrise the positive roots of the dual \smash{$G_2^{(1)}$} of \smash{$D_4^{(3)}$}, with Cartan matrix shown in Figure~\ref{fig:CMs}. The basic imaginary root is given by $\delta = \alpha_0 + 2\alpha_1 + 3\alpha_2$. The set of real roots is given by
\[
\Phi_{\mathrm{Re}} = \{\alpha + r\delta \mid \alpha\in\Phi,\, r\in \Z\},
\]
and the set of imaginary roots is given by
\[
\Phi_{\mathrm{Im}} = \{k\delta \mid k\in \Z \backslash 0\}.
\]
Therefore, we can parametrise the set of positive roots of $G_2^{(1)}$ as follows. The short, real positive roots are given by
\begin{align*}
\Phi_{\mathrm{Re},s}^{+} &=
\left\{
\begin{array}{@{}ll@{}}
 \alpha_2 + k(\alpha_0 + 2\alpha_1 + 3\alpha_2) & \text{for} \ k \geq 0 \\
 \alpha_1 + \alpha_2 + k(\alpha_0 + 2\alpha_1 + 3\alpha_2) & \text{for} \  k \geq 0 \\
 \alpha_1 + 2\alpha_2 + k(\alpha_0 + 2\alpha_1 + 3\alpha_2) & \text{for} \  k \geq 0 \\
 -\alpha_2 + k(\alpha_0 + 2\alpha_1 + 3\alpha_2) & \text{for} \  k \geq 1 \\
 -\alpha_1 - \alpha_2 + k(\alpha_0 + 2\alpha_1 + 3\alpha_2) & \text{for} \  k \geq 1 \\
 -\alpha_1 - 2\alpha_2 + k(\alpha_0 + 2\alpha_1 + 3\alpha_2) & \text{for} \  k \geq 1
\end{array}
\right\} \\
&=
\left\{
\begin{array}{@{}ll@{}}
 k\alpha_0 + 2k\alpha_1 + (3k+1)\alpha_2 & \text{for} \  k \geq 0 \\
 k\alpha_0 + (2k+1)\alpha_1 + (3k+1)\alpha_2 & \text{for} \  k \geq 0 \\
 k\alpha_0 + (2k+1)\alpha_1 + (3k+2)\alpha_2 & \text{for} \  k \geq 0 \\
 k\alpha_0 + 2k\alpha_1 + (3k-1)\alpha_2 & \text{for} \  k \geq 1 \\
 k\alpha_0 + (2k-1)\alpha_1 + (3k-1)\alpha_2 & \text{for} \  k \geq 1 \\
 k\alpha_0 + (2k-1)\alpha_1 + (3k-2)\alpha_2 & \text{for} \  k \geq 1
\end{array}
\right\},
\end{align*}
the long, real positive roots are given by
\begin{align*}
\Phi_{\mathrm{Re},l}^{+} &=
\left\{
\begin{array}{@{}ll@{}}
 \alpha_1 + k(\alpha_0 + 2\alpha_1 + 3\alpha_2) & \text{for} \  k \geq 0 \\
 \alpha_1 + 3\alpha_2 + k(\alpha_0 + 2\alpha_1 + 3\alpha_2) & \text{for} \  k \geq 0 \\
 2\alpha_1 + 3\alpha_2 + k(\alpha_0 + 2\alpha_1 + 3\alpha_2) & \text{for} \  k \geq 0 \\
 -\alpha_1 + k(\alpha_0 + 2\alpha_1 + 3\alpha_2) & \text{for} \  k \geq 1 \\
 -\alpha_1 - 3\alpha_2 + k(\alpha_0 + 2\alpha_1 + 3\alpha_2) & \text{for} \  k \geq 1 \\
 -2\alpha_1 - 3\alpha_2 + k(\alpha_0 + 2\alpha_1 + 3\alpha_2) & \text{for} \  k \geq 1
\end{array}
\right\} \\
&=
\left\{
\begin{array}{@{}ll@{}}
 k\alpha_0 + (2k+1)\alpha_1 + 3k\alpha_2 & \text{for} \  k \geq 0 \\
 k\alpha_0 + (2k+1)\alpha_1 + (3k+3)\alpha_2 & \text{for} \  k \geq 0 \\
 k\alpha_0 + (2k+2)\alpha_1 + (3k+3)\alpha_2 & \text{for} \  k \geq 0 \\
 k\alpha_0 + (2k-1)\alpha_1 + 3k\alpha_2 & \text{for} \  k \geq 1 \\
 k\alpha_0 + (2k-1)\alpha_1 + (3k-3)\alpha_2 & \text{for} \  k \geq 1 \\
 k\alpha_0 + (2k-2)\alpha_1 + (3k-3)\alpha_2 & \text{for} \  k \geq 1
\end{array}
\right\},
\end{align*}
and the positive imaginary roots are given by
\begin{align*}\Phi_{\mathrm{Im}}^{+} =
\{
k(\alpha_0 + 2\alpha_1 + 3\alpha_2) \ \text{for} \ k \geq 1
\}.
\end{align*}
For $\alpha\in\Phi_{\mathrm{Re}}$, we have $m(\alpha)=1$. For $\alpha\in\Phi_{\mathrm{Im}}$ in type \smash{$G_2^{(1)}$}, we have $m(\alpha)=2$.

We need to compute specialisations of the product
\begin{align*}
D(\Phi^\vee):={}&\prod_{\alpha \in \Phi^+} (1 - {\rm e}^{-\alpha})^{m(\alpha)} \\
={}&
\prod_{\alpha \in \Phi_{\mathrm{Re},s}^+} (1 - {\rm e}^{-\alpha})^{m(\alpha)}
\prod_{\alpha \in \Phi_{\mathrm{Re},l}^+} (1 - {\rm e}^{-\alpha})^{m(\alpha)}
\prod_{\alpha \in \Phi_{\mathrm{Im}}^+} (1 - {\rm e}^{-\alpha})^{m(\alpha)},
\end{align*}
which takes the following form. For the real short roots, this product is
\begin{align*}
&\prod_{\alpha \in \Phi_{\mathrm{Re},s}^+}(1 - {\rm e}^{-\alpha})^{m(\alpha)}
=\prod_{k=0}^{\infty} \bigl(1 - {\rm e}^{-k\alpha_0} {\rm e}^{-2k\alpha_1} {\rm e}^{-(3k+1)\alpha_2}\bigr)
\\
&\qquad{} \times \prod_{k=0}^{\infty} \bigl(1 - {\rm e}^{-k\alpha_0} {\rm e}^{-(2k+1)\alpha_1} {\rm e}^{-(3k+1)\alpha_2}\bigr)
\prod_{k=1}^{\infty} \bigl(1 - {\rm e}^{-k\alpha_0} {\rm e}^{-2k\alpha_1} {\rm e}^{-(3k-1)\alpha_2}\bigr) \\
&\qquad{} \times \prod_{k=0}^{\infty} \bigl(1 - {\rm e}^{-k\alpha_0} {\rm e}^{-(2k+1)\alpha_1} {\rm e}^{-(3k+2)\alpha_2}\bigr)\prod_{k=1}^{\infty} \bigl(1 - {\rm e}^{-k\alpha_0} {\rm e}^{-(2k-1)\alpha_1} {\rm e}^{-(3k-1)\alpha_2}\bigr)\\
&\qquad{}\times \prod_{k=1}^{\infty} \bigl(1 - {\rm e}^{-k\alpha_0} {\rm e}^{-(2k-1)\alpha_1} {\rm e}^{-(3k-2)\alpha_2}\bigr)
\end{align*}
for real long roots, we have
\begin{align*}
&\prod_{\alpha \in \Phi_{\mathrm{Re},l}^+} (1 - {\rm e}^{-\alpha})^{m(\alpha)}
=\prod_{k=0}^{\infty} \bigl(1 - {\rm e}^{-k\alpha_0} {\rm e}^{-(2k+1)\alpha_1} {\rm e}^{-3k\alpha_2}\bigr) \\
&\qquad \times \prod_{k=0}^{\infty} \bigl(1 - {\rm e}^{-k\alpha_0} {\rm e}^{-(2k+1)\alpha_1} {\rm e}^{-(3k+3)\alpha_2}\bigr)\prod_{k=0}^{\infty} \bigl(1 - {\rm e}^{-k\alpha_0} {\rm e}^{-(2k+2)\alpha_1} {\rm e}^{-(3k+3)\alpha_2}\bigr)\\
&\qquad{} \times \prod_{k=1}^{\infty} \bigl(1 - {\rm e}^{-k\alpha_0} {\rm e}^{-(2k-1)\alpha_1} {\rm e}^{-3k\alpha_2}\bigr)\prod_{k=1}^{\infty} \bigl(1 - {\rm e}^{-k\alpha_0} {\rm e}^{-(2k-1)\alpha_1} {\rm e}^{-(3k-3)\alpha_2}\bigr) \\
&\qquad{} \times \prod_{k=1}^{\infty} \bigl(1 - {\rm e}^{-k\alpha_0} {\rm e}^{-(2k-2)\alpha_1} {\rm e}^{-(3k-3)\alpha_2}\bigr),
\end{align*}
and for the imaginary roots, we have
\begin{align*}
\prod_{\alpha \in \Phi_{\mathrm{Im}}^+} (1 - {\rm e}^{-\alpha} )^{m(\alpha)}
&=
\prod_{k=1}^{\infty} \bigl(1 - {\rm e}^{-k\alpha_0} {\rm e}^{-2k\alpha_1} {\rm e}^{-3k\alpha_2}\bigr)^2.
\end{align*}

Thus, the principal specialisation of $D\bigl(\Phi^\vee\bigr)$ simplifies to the following product of infinite $q$-Pochhammer symbols:
\begin{align*}
F_1\!\! \prod_{\alpha \in \Phi^+}\! (1 - {\rm e}^{-\alpha} )^{m(\alpha)}\!
={}&
\bigl(q;q^6\bigr)_{\infty} \bigl(q^2;q^6\bigr)_{\infty} \bigl(q^3;q^6\bigr)_{\infty} \bigl(q^5;q^6\bigr)_{\infty} \bigl(q^4;q^6\bigr)_{\infty} \bigl(q^3;q^6\bigr)_{\infty} \\[-2.3ex]
&{}{\times}  \bigl(q;q^6\bigr)_{\infty} \bigl(q^4;q^6\bigr)_{\infty} \bigl(q^5;q^6\bigr)_{\infty} \bigl(q^5;q^6\bigr)_{\infty} \bigl(q^2;q^6\bigr)_{\infty} \bigl(q;q^6\bigr)_{\infty} \bigl(q^6;q^6\bigr)^2 \\
={}&
\bigl(q;q\bigr)_{\infty}^2 \bigl(q^5;q^6\bigr)_{\infty} \bigl(q;q^6\bigr)_{\infty},
\end{align*}
and the $(2,1,1)$-specialisation of $D\bigl(\Phi^\vee\bigr)$ simplifies to
\begin{align*}
F_{(2,1,1)}\!\! \prod_{\alpha \in \Phi^+}\! (1 - {\rm e}^{-\alpha})^{m(\alpha)}
\!={}&
\bigl(q;q^7\bigr)_{\infty} \bigl(q^2;q^7\bigr)_{\infty} \bigl(q^3;q^7\bigr)_{\infty} \bigl(q^6;q^7\bigr)_{\infty} \bigl(q^5;q^7\bigr)_{\infty} \bigl(q^4;q^7\bigr)_{\infty} \bigl(q;q^7\bigr)_{\infty} \\[-2.3ex]
&{}{\times}  \bigl(q^4;q^7\bigr)_{\infty} \bigl(q^5;q^7\bigr)_{\infty} \bigl(q^6;q^7\bigr)_{\infty} \bigl(q^3;q^7\bigr)_{\infty} \bigl(q^2;q^7\bigr)_{\infty} \bigl(q^7;q^7\bigr)^2 \\
={}&
\bigl(q;q\bigr)_{\infty}^2.
\end{align*}

Hence, the principal specialisation of the highest weight $\Lambda_0$-character is given by
\begin{align}\label{eq:cong_cond}
 F_1\bigl({\rm e}^{-\Lambda_0}\mathrm{ch} (L(\Lambda_0) )\bigr) = \frac{(q;q)_{\infty}^2}{(q;q)_{\infty}^2\bigl(q^5;q^6\bigr)_{\infty}\bigl(q;q^6\bigr)_{\infty}}= \frac{1}{\bigl(q,q^5;q^6\bigr)_{\infty}},
\end{align}
the generating function of partitions with parts congruent to $1$ or $5$ modulo $6$.

\section{Difference conditions}\label{section:difference}

In this section, we express the character of the irreducible representation of highest weight~$\Lambda_0$~of the twisted affine Kac--Moody algebra $\mathfrak{g}$ of type \smash{$D_4^{(3)}$}, using a perfect crystal of type \smash{$D_4^{(3)}$} and the combinatorial framework of Dousse--Konan \cite{DK19b,DK22}.

\begin{Theorem}[{Dousse--Konan \cite[Theorem~1.9]{DK19b}}]\label{thm:KMN2coloured}
Let $\mathcal{B}$ be a perfect crystal of level $\ell$ for an~affine Kac--Moody algebra $\mathfrak{g}$ with vertices $v_1,v_2,\ldots,v_n$. Let $\lambda \in P^+_\ell$ with a constant ground state path $p_\lambda = \cdots \otimes v_j \otimes v_j \otimes v_j$ for some $j\in\{1,\ldots,n\}$. Let $H$ be the energy matrix of~$\mathcal{B}$ determined by the initial condition $H(v_j\otimes v_j)=0$, with rows and columns indexed by colours~${c_1,c_2,\ldots,c_n}$. To a grounded partition $\pi=(\pi_1,\pi_2,\ldots,\pi_\ell)$ where the colour of the part $\pi_j$ is $c_{i_j}$ for $j\in\{1,\ldots,\ell\}$, we associate the monomial $C(\pi):= c_{i_1}c_{i_2} \cdots  c_{i_\ell}$.
Setting \smash{$t = {\rm e}^{-\delta / a_0}$}, and setting \smash{$c_v = {\rm e}^{\overline{\mathrm{wt}}\, v}$} for all~${v \in \mathcal{B}}$, we have
\begin{align}
&\sum_{\pi \in \mathcal{P}_{H}^{c_j}} C(\pi) t^{|\pi|}= {\rm e}^{-\lambda} \operatorname{ch}L(\lambda),\label{eq:grounded_char}\\
&\sum_{\pi \in \mathcal{P}_{H}^{c_j,\geq}} C(\pi) t^{|\pi|}= \frac{{\rm e}^{-\lambda} \operatorname{ch}L(\lambda)}{(t;t)_\infty}\label{eq:geq_grounded_char}.
\end{align}
\end{Theorem}

We will use this combinatorial character formula on the perfect crystal corresponding to the level 1 irreducible representation of highest weight $\Lambda_0$ of $\mathfrak{g}$. This crystal was explicitly constructed by Kashiwara, Misra, Okado and Yamada~\cite{KMOY07}, shown in Figure~\ref{fig:G2perfectcrystal}, where the $1$\nobreakdash-arrows are coloured blue, the $2$-arrows are coloured red, and the $0$-arrows are coloured green.

\begin{figure}[!t]
 \centering
 \begin{tikzpicture}[scale=0.75, every node/.style={scale=1.2}]
 \node (a1) at (0,0) {$v_a$};
 \node (a2) at (1.5,0) {$v_b$};
 \node (a3) at (3,0) {$v_c$};
 \node (a0) at (4.5,0) {$v_d$};
 \node (b3) at (6,0) {$v_e$};
 \node (b2) at (7.5,0) {$v_f$};
 \node (b1) at (9,0) {$v_g$};
 \node (phi) at (4.5,-1.5) {$v_h$};

 \draw[->,blue,thick] (a1) -- (a2);
 \draw[->,blue,thick] (a3) -- (a0);
 \draw[->,blue,thick] (b2) -- (b1);
 \draw[->,blue,thick] (a0) -- (b3);

 \draw[->,red,thick] (a2) -- (a3);
 \draw[->,red,thick] (b3) -- (b2);

 \draw[->,green!60!black,thick] (b1) to [bend left=15] (phi);
 \draw[->,green!60!black,thick] (phi) to [bend left=15] (a1);
 \draw[->,green!60!black,thick] (b3) to [bend right=30] (a2);
 \draw[->,green!60!black,thick] (b2) to [bend right=30] (a3);
 \end{tikzpicture}
 \caption{Perfect crystal $\mathcal{B}$ of level 1 of type $D_4^{(3)}$.} \label{fig:G2perfectcrystal}
\end{figure}

We compute the levels of vertices $v\in\mathcal{B}$, using that $c= h_0 + 2h_1 + 3h_2$ in $\mathfrak{g}$:
\begin{align*}
(\varepsilon(v))(c) &= (\varepsilon(v))(h_0 + 2h_1 + 3h_2) = \varepsilon_0(v) + 2 \varepsilon_1(v) + 3 \varepsilon_2(v).
\end{align*}
Now we can read off the values of $\varphi(v)$ and $\varepsilon(v)$ for each vertex $v\in\mathcal{B}$, compute the weights $\overline{\mathrm{wt}}(v) = \varphi(v) - \varepsilon(v)$, and the levels. We summarise these values in Figure~\ref{fig:wts}.
{\small\begin{figure}[!t]
 \centering
 \begin{tabular}{|c|c|c|c|c|}
\hline
$v$ & \text{Level} & $\varphi(v)$ & $\varepsilon(v)$ & \text{Weight} \\
\hline
$v_a$ & 2 & $\Lambda_1$ & $2\Lambda_0$ & $-2\Lambda_0 + \Lambda_1$ \\
$v_b$ & 3 & $\Lambda_2$ & $\Lambda_0 + \Lambda_1$ & $-\Lambda_0 - \Lambda_1 + \Lambda_2$ \\
$v_c$ & 4 & $2\Lambda_1$ & $\Lambda_0 + \Lambda_2$ & $-\Lambda_0 + 2\Lambda_1 - \Lambda_2$ \\
$v_d$ & 2 & $\Lambda_1$ & $\Lambda_1$ & $0$ \\
$v_e$ & 4 & $\Lambda_0 + \Lambda_2$ & $2\Lambda_1$ & $\Lambda_0 - 2\Lambda_1 + \Lambda_2$ \\
$v_f$ & 3 & $\Lambda_0 + \Lambda_1$ & $\Lambda_2$ & $\Lambda_0 + \Lambda_1 - \Lambda_2$ \\
$v_g$ & 2 & $2\Lambda_0$ & $\Lambda_1$ & $2\Lambda_0 - \Lambda_1$ \\
$v_h$ & 1 & $\Lambda_0$ & $\Lambda_0$ & $0$ \\
\hline
\end{tabular}
\caption{Weights and levels of the vertices of $\mathcal{B}$.}
 \label{fig:wts}
\end{figure}}
For this crystal to be perfect, we must first verify that there are unique elements $v^\lambda$ and $v_\lambda$ in $\mathcal{B}$ such that $\varepsilon(v^\lambda) = \Lambda_0$ and $\varphi(v_\lambda) = \Lambda_0$. Indeed, we have $v^\lambda = v_\lambda = v_h$, resulting in a constant ground state path of level~1, $\cdots\otimes v_h\otimes v_h\otimes v_h$. Furthermore, the restriction of the maps $\varepsilon$ and $\varphi$ to $P_\ell^+$ is bijective, as there is only one level-1 vertex in $\mathcal{B}$. Next, we need to compute $\mathcal{B} \otimes \mathcal{B}$ to verify that it is connected. The explicit description of $\mathcal{B} \otimes \mathcal{B}$ will also be necessary to compute the energy matrix. By Definition \ref{def:crystaltensor}, the arrows in $\mathcal{B} \otimes \mathcal{B}$ can be computed independently. The finite part, which we obtain from $\mathcal{B} \otimes \mathcal{B}$ by deleting the $0$-arrows, has six connected components: two copies of the trivial representation of $\mathfrak{g}_\mathrm{fin} = G_2$, three copies of the finite part of $\mathcal{B}$, one copy of the $14$-dimensional regular representation of $\mathfrak{g}_\mathrm{fin}$ and a $27$-dimensional irreducible representation. On the other hand, the crystal $\mathcal{B}\otimes\mathcal{B}$, including the 0-arrows (in green) is connected. This is shown in Figure~\ref{fig:G2tensor}. Using Definition \ref{def:energy}, we can directly read off the energy matrix $H$ from the crystal $\mathcal{B} \otimes \mathcal{B}$, shown in Figure~\ref{fig:G2energyfn}.

\begin{figure}[!t]
 \centering
 \begin{tikzpicture}[scale=0.75, every node/.style={scale=0.9}]
 \node (n11) at (0,0) {$v_a \otimes v_a$};
 \node (n12) at (2.5,0) {$v_a \otimes v_b$};
 \node (n13) at (5,0) {$v_a \otimes v_c$};
 \node (n14) at (7.5,0) {$v_a \otimes v_d$};
 \node (n15) at (10,0) {$v_a \otimes v_e$};
 \node (n16) at (12.5,0) {$v_a \otimes v_f$};
 \node (n17) at (15,0) {$v_a \otimes v_g$};
 \node (n18) at (17.5,0) {$v_a \otimes v_h$};

 \node (n21) at (0,-2) {$v_b \otimes v_a$};
 \node (n22) at (2.5,-2) {$v_b \otimes v_b$};
 \node (n23) at (5,-2) {$v_b \otimes v_c$};
 \node (n24) at (7.5,-2) {$v_b \otimes v_d$};
 \node (n25) at (10,-2) {$v_b \otimes v_e$};
 \node (n26) at (12.5,-2) {$v_b \otimes v_f$};
 \node (n27) at (15,-2) {$v_b \otimes v_g$};
 \node (n28) at (17.5,-2) {$v_b \otimes v_h$};

 \node (n31) at (0,-4) {$v_c \otimes v_a$};
 \node (n32) at (2.5,-4) {$v_c \otimes v_b$};
 \node (n33) at (5,-4) {$v_c \otimes v_c$};
 \node (n34) at (7.5,-4) {$v_c \otimes v_d$};
 \node (n35) at (10,-4) {$v_c \otimes v_e$};
 \node (n36) at (12.5,-4) {$v_c \otimes v_f$};
 \node (n37) at (15,-4) {$v_c \otimes v_g$};
 \node (n38) at (17.5,-4) {$v_c \otimes v_h$};

 \node (n41) at (0,-6) {$v_d \otimes v_a$};
 \node (n42) at (2.5,-6) {$v_d \otimes v_b$};
 \node (n43) at (5,-6) {$v_d \otimes v_c$};
 \node (n44) at (7.5,-6) {$v_d \otimes v_d$};
 \node (n45) at (10,-6) {$v_d \otimes v_e$};
 \node (n46) at (12.5,-6) {$v_d \otimes v_f$};
 \node (n47) at (15,-6) {$v_d \otimes v_g$};
 \node (n48) at (17.5,-6) {$v_d \otimes v_h$};

 \node (n51) at (0,-8) {$v_e \otimes v_a$};
 \node (n52) at (2.5,-8) {$v_e \otimes v_b$};
 \node (n53) at (5,-8) {$v_e \otimes v_c$};
 \node (n54) at (7.5,-8) {$v_e \otimes v_d$};
 \node (n55) at (10,-8) {$v_e \otimes v_e$};
 \node (n56) at (12.5,-8) {$v_e \otimes v_f$};
 \node (n57) at (15,-8) {$v_e \otimes v_g$};
 \node (n58) at (17.5,-8) {$v_e \otimes v_h$};

 \node (n61) at (0,-10) {$v_f \otimes v_a$};
 \node (n62) at (2.5,-10) {$v_f \otimes v_b$};
 \node (n63) at (5,-10) {$v_f \otimes v_c$};
 \node (n64) at (7.5,-10) {$v_f \otimes v_d$};
 \node (n65) at (10,-10) {$v_f \otimes v_e$};
 \node (n66) at (12.5,-10) {$v_f \otimes v_f$};
 \node (n67) at (15,-10) {$v_f \otimes v_g$};
 \node (n68) at (17.5,-10) {$v_f \otimes v_h$};

 \node (n71) at (0,-12) {$v_g \otimes v_a$};
 \node (n72) at (2.5,-12) {$v_g \otimes v_b$};
 \node (n73) at (5,-12) {$v_g \otimes v_c$};
 \node (n74) at (7.5,-12) {$v_g \otimes v_d$};
 \node (n75) at (10,-12) {$v_g \otimes v_e$};
 \node (n76) at (12.5,-12) {$v_g \otimes v_f$};
 \node (n77) at (15,-12) {$v_g \otimes v_g$};
 \node (n78) at (17.5,-12) {$v_g \otimes v_h$};

 \node (n81) at (0,-14) {$v_h \otimes v_a$};
 \node (n82) at (2.5,-14) {$v_h \otimes v_b$};
 \node (n83) at (5,-14) {$v_h \otimes v_c$};
 \node (n84) at (7.5,-14) {$v_h \otimes v_d$};
 \node (n85) at (10,-14) {$v_h \otimes v_e$};
 \node (n86) at (12.5,-14) {$v_h \otimes v_f$};
 \node (n87) at (15,-14) {$v_h \otimes v_g$};
 \node (n88) at (17.5,-14) {$v_h \otimes v_h$};

 \node (n91) at (0,-16) {};
 \node (n92) at (2.5,-16) {};
 \node (n93) at (5,-16) {};
 \node (n94) at (7.5,-16) {};
 \node (n95) at (10,-16) {};
 \node (n96) at (12.5,-16) {};
 \node (n97) at (15,-16) {};
 \node (n98) at (17.5,-16) {};

 \node (n19) at (20,0) {};
 \node (n29) at (20,-2) {};
 \node (n39) at (20,-4) {};
 \node (n49) at (20,-6) {};
 \node (n59) at (20,-8) {};
 \node (n69) at (20,-10) {};
 \node (n79) at (20,-12) {};
 \node (n89) at (20,-14) {};
 \node (n99) at (20,-16) {};

 \draw[->,red,thick] (n12) -- (n13);
 \draw[->,red,thick] (n32) -- (n33);
 \draw[->,red,thick] (n42) -- (n43);
 \draw[->,red,thick] (n62) -- (n63);
 \draw[->,red,thick] (n72) -- (n73);
 \draw[->,red,thick] (n82) -- (n83);

 \draw[->,red,thick] (n15) -- (n16);
 \draw[->,red,thick] (n35) -- (n36);
 \draw[->,red,thick] (n45) -- (n46);
 \draw[->,red,thick] (n65) -- (n66);
 \draw[->,red,thick] (n75) -- (n76);
 \draw[->,red,thick] (n85) -- (n86);

 \draw[->,red,thick] (n21) -- (n31);
 \draw[->,red,thick] (n22) -- (n32);
 \draw[->,red,thick] (n24) -- (n34);
 \draw[->,red,thick] (n25) -- (n35);
 \draw[->,red,thick] (n27) -- (n37);
 \draw[->,red,thick] (n28) -- (n38);

 \draw[->,red,thick] (n51) -- (n61);
 \draw[->,red,thick] (n52) -- (n62);
 \draw[->,red,thick] (n54) -- (n64);
 \draw[->,red,thick] (n55) -- (n65);
 \draw[->,red,thick] (n57) -- (n67);
 \draw[->,red,thick] (n58) -- (n68);

 \draw[->,blue,thick] (n14) -- (n15);

 \draw[->,blue,thick] (n21) -- (n22);
 \draw[->,blue,thick] (n23) -- (n24);
 \draw[->,blue,thick] (n24) -- (n25);
 \draw[->,blue,thick] (n26) -- (n27);

 \draw[->,blue,thick] (n44) -- (n45);

 \draw[->,blue,thick] (n51) -- (n52);
 \draw[->,blue,thick] (n53) -- (n54);
 \draw[->,blue,thick] (n54) -- (n55);
 \draw[->,blue,thick] (n56) -- (n57);

 \draw[->,blue,thick] (n64) -- (n65);

 \draw[->,blue,thick] (n71) -- (n72);
 \draw[->,blue,thick] (n73) -- (n74);
 \draw[->,blue,thick] (n74) -- (n75);
 \draw[->,blue,thick] (n76) -- (n77);

 \draw[->,blue,thick] (n81) -- (n82);
 \draw[->,blue,thick] (n83) -- (n84);
 \draw[->,blue,thick] (n84) -- (n85);
 \draw[->,blue,thick] (n86) -- (n87);

 \draw[->,blue,thick] (n11) -- (n21);
 \draw[->,blue,thick] (n31) -- (n41);
 \draw[->,blue,thick] (n41) -- (n51);
 \draw[->,blue,thick] (n61) -- (n71);

 \draw[->,blue,thick] (n32) -- (n42);

 \draw[->,blue,thick] (n13) -- (n23);
 \draw[->,blue,thick] (n33) -- (n43);
 \draw[->,blue,thick] (n43) -- (n53);
 \draw[->,blue,thick] (n63) -- (n73);

 \draw[->,blue,thick] (n34) -- (n44);

 \draw[->,blue,thick] (n16) -- (n26);
 \draw[->,blue,thick] (n36) -- (n46);
 \draw[->,blue,thick] (n46) -- (n56);
 \draw[->,blue,thick] (n66) -- (n76);

 \draw[->,blue,thick] (n37) -- (n47);

 \draw[->,blue,thick] (n18) -- (n28);
 \draw[->,blue,thick] (n38) -- (n48);
 \draw[->,blue,thick] (n48) -- (n58);
 \draw[->,blue,thick] (n68) -- (n78);

 \draw[->,green!60!black,thick] (n16) to [bend right=15] (n13);
 \draw[->,green!60!black,thick] (n15) to [bend right=15] (n12);

 \draw[->,green!60!black,thick] (n26) to [bend right=15] (n23);
 \draw[->,green!60!black,thick] (n25) to [bend right=15] (n22);

 \draw[->,green!60!black,thick] (n36) to [bend right=15] (n33);
 \draw[->,green!60!black,thick] (n35) to [bend right=15] (n32);

 \draw[->,green!60!black,thick] (n46) to [bend right=15] (n43);
 \draw[->,green!60!black,thick] (n45) to [bend right=15] (n42);

 \draw[->,green!60!black,thick] (n18) -- (n19);
 \draw[->,green!60!black,thick] (n28) -- (n29);
 \draw[->,green!60!black,thick] (n38) -- (n39);
 \draw[->,green!60!black,thick] (n48) -- (n49);
 \draw[->,green!60!black,thick] (n58) -- (n59);
 \draw[->,green!60!black,thick] (n68) -- (n69);
 \draw[->,green!60!black,thick] (n88) -- (n89);

 \draw[->,green!60!black,thick] (n64) to [bend left=27] (n34);
 \draw[->,green!60!black,thick] (n54) to [bend left=27] (n24);

 \draw[->,green!60!black,thick] (n65) to [bend left=27] (n35);
 \draw[->,green!60!black,thick] (n55) to [bend left=27] (n25);

 \draw[->,green!60!black,thick] (n66) to [bend left=27] (n36);
 \draw[->,green!60!black,thick] (n56) to [bend left=27] (n26);

 \draw[->,green!60!black,thick] (n67) to [bend left=27] (n37);
 \draw[->,green!60!black,thick] (n57) to [bend left=27] (n27);

 \draw[->,green!60!black,thick] (n84) -- (n94);
 \draw[->,green!60!black,thick] (n85) -- (n95);
 \draw[->,green!60!black,thick] (n86) -- (n96);
 \draw[->,green!60!black!60!black,thick] (n87) -- (n97);
\end{tikzpicture}

\caption{the crystal $\mathcal{B} \otimes \mathcal{B}$.} \label{fig:G2tensor}
\end{figure}

\begin{figure}[t]
\centering
$
\begin{array}{c|cccccccc}
& v_a & v_b & v_c & v_d & v_e & v_f & v_g & v_h \\
\hline
v_a & {\color{black}2} & {\color{black}1} & {\color{black}1} & {\color{black}0} & {\color{black}0} & {\color{black}0} & {\color{black}0} & {\color{black}1} \\
v_b & {\color{black}2} & {\color{black}2} & {\color{black}1} & {\color{black}1} & {\color{black}1} & {\color{black}0} & {\color{black}0} & {\color{black}1} \\
v_c & {\color{black}2} & {\color{black}2} & {\color{black}2} & {\color{black}1} & {\color{black}1} & {\color{black}1} & {\color{black}0} & {\color{black}1} \\
v_d & {\color{black}2} & {\color{black}2} & {\color{black}2} & {\color{black}1} & {\color{black}1} & {\color{black}1} & {\color{black}0} & {\color{black}1} \\
v_e & {\color{black}2} & {\color{black}2} & {\color{black}2} & {\color{black}2} & {\color{black}2} & {\color{black}1} & {\color{black}1} & {\color{black}1} \\
v_f & {\color{black}2} & {\color{black}2} & {\color{black}2} & {\color{black}2} & {\color{black}2} & {\color{black}2} & {\color{black}1} & {\color{black}1} \\
v_g & {\color{black}2} & {\color{black}2} & {\color{black}2} & {\color{black}2} & {\color{black}2} & {\color{black}2} & {\color{black}2} & {\color{black}1} \\
v_h & {\color{black}1} & {\color{black}1} & {\color{black}1} & {\color{black}1} & {\color{black}1} & {\color{black}1} & {\color{black}1} & {\color{black}0} \\
\end{array}
$
\caption{Energy matrix $H$ of $\mathcal{B}$.}
\label{fig:G2energyfn}
\end{figure}

Recall from Figure~\ref{fig:wts} that the two 0-weight vertices are $v_d$ and $v_h$. To express the weights~$\overline{\mathrm{wt}}(v)$ for all $v\in\mathcal{B}$ in terms of the simple roots $\alpha_1$ and $\alpha_2$, we simply need to follow the arrows in the crystal. For $i \in \{0,1,2\}$, an arrow labelled $i$ decreases the weight by $\alpha_i$. Note that $\delta = \alpha_0 + 2\alpha_1 + \alpha_2 = 0$ in the derived algebra $\mathfrak{g}'$, which is consistent with the cycles in the crystal graph (Figure~\ref{fig:G2perfectcrystal}). The assumption $\lambda(d)=0$, together with condition~(iii) of Definition~\ref{def:perfectcrystal}, implies that for each $v\in\mathcal{B}$ the coefficient of $\alpha_0$ in $\overline{\mathrm{wt}}(v)$ is zero. Starting from vertex $v_d$, which has weight $0$, we have: $\overline{\mathrm{wt}}(v_e)= -\alpha_1$, $\overline{\mathrm{wt}}(v_f)= -\alpha_1-\alpha_2$, $\overline{\mathrm{wt}}(v_g)= -2\alpha_1-\alpha_2$, $\overline{\mathrm{wt}}(v_h)= 0$, $\overline{\mathrm{wt}}(v_a)= 2\alpha_1+\alpha_2$, $\overline{\mathrm{wt}}(v_b)= \alpha_1+\alpha_2$ and $\overline{\mathrm{wt}}(v_c)= \alpha_1$.

By Theorem~\ref{thm:KMN2coloured}, equation~\eqref{eq:grounded_char},
the character can be written as
\begin{align*}
\sum_{\pi \in \mathcal{P}_H^{h}} C(\pi) t^{|\pi|} = {\rm e}^{-\Lambda_0} \operatorname{ch} L(\Lambda_0),
\end{align*}
where $t = {\rm e}^{-\delta}$, \smash{$c_v = {\rm e}^{\overline{\mathrm{wt}}\, v}$} for all $v \in \mathcal{B}$, and $\mathcal{P}_H^{h}$ is the set of grounded partitions with ground~$0_{h}$ with respect to the energy function shown in Figure~\ref{fig:G2energyfn}. It can be seen from the generalised~Cartan matrix in Figure~\ref{fig:CMs} that $a_0=1$ in type \smash{$D_4^{(3)}$}.

The ground state path corresponding to the highest weight $\Lambda_0$ is
\[
p_{\Lambda_0}=\cdots\otimes v_h\otimes v_h\otimes v_h.
\]

Three examples of paths in $P(\Lambda_0)$ are
\begin{equation}\label{eq:paths}
\begin{gathered}
\cdots\otimes v_h\otimes v_h\otimes v_h\otimes v_f\otimes v_e\otimes v_d\otimes v_c\otimes v_b\otimes v_h,\\
\cdots\otimes v_h\otimes v_h\otimes v_g\otimes v_f\otimes v_f\otimes v_e\otimes v_d\otimes v_b\otimes v_h,\\
\cdots\otimes v_h\otimes v_h\otimes v_h\otimes v_g\otimes v_e\otimes v_e\otimes v_c\otimes v_b\otimes v_a.
\end{gathered}
\end{equation}
To avoid confusion, we introduce the variables $a$, $b$, $c$, $d$, $e$, $f$, $g$, $h$ for the colours corresponding to the crystal vertices $v_a$, $v_b$, $v_c$, $v_d$, $v_e$, $v_f$, $v_g$, $v_h$, respectively. The grounded partitions~corresponding to the paths in~\eqref{eq:paths} are
\begin{align*}
&0_h\,1_f\,2_e\,3_d\,4_c\,5_b\,6_h,\\
&0_h\,1_g\,2_f\,4_f\,5_e\,6_d\,7_b\,8_h,\\
&0_h\,1_g\,2_e\,4_e\,5_c\,6_b\,7_a.
\end{align*}

In the principal specialisation, for each simple root $\alpha_i$, ${\rm e}^{-\alpha_i}$ is mapped to $q$. In particular, we have the dilation $t = q^4$. The weights, expressed in terms of $\alpha_1$ and $\alpha_2$, along with their corresponding colour specialisations, are listed in Figure~\ref{fig:specialisations}.
\begin{figure}[!t]
\centering
 \begin{tabular}{|c|c|}
 \hline
 & $\overline{wt}(b)$ \\
 \hline
 $v_a$ & $ 2\alpha_1+\alpha_2$ \\
 $v_b$ & $ \alpha_1+\alpha_2$ \\
 $v_c$ & $\alpha_1$ \\
 $v_d$ & $0$ \\
 $v_e$ & $-\alpha_1$ \\
 $v_f$ & $-\alpha_1-\alpha_2$ \\
 $v_g$ & $-2\alpha_1-\alpha_2$ \\
 $v_h$ & $0$ \\
 \hline
 \end{tabular}
 \qquad
 {$\longrightarrow$}
 \qquad
 $ \begin{aligned}
 &\underline{\text{specialisations}}
 \\
 &d, h \mapsto 1 \\
 &a \mapsto q^{-3} \\
 &b \mapsto q^{-2} \\
 &g \mapsto q^3 \\
 &f \mapsto q^2 \\
 &e \mapsto q \\
 &c \mapsto q^{-1}
 \end{aligned}$

 \caption{Specialisations of the weights.}
\label{fig:specialisations}
\end{figure}

We aim to find a square matrix $M$ with rows and columns indexed by $a$, $b$, $c$, $d$, $e$, $f$, $g$, $h$ such that the specialisation of the sum in (\ref{thm:KMN2coloured}) becomes the following generating function:
\begin{align*}
F_1\biggl(\sum_{\pi \in \mathcal{P}_{H}^h} C(\pi) t^{|\pi|}\biggr) = \sum_{\pi \in \mathcal{P}_{M}^h}q^{|p|}.
\end{align*}

To this end, we need to alter the matrix according to the specialisations above. For example, the entry in column $b$ and row $a$ in $H$ is $1$, which means that if a part with colour $a$ follows a part with colour $b$, then their difference is $1$. The relevant specialisations are $t = q^4$, $a = q^{-3}$, and $b = q^{-2}$, and the corresponding value in the new matrix $M$ should be: $1 \cdot 4 + (-2) - (-3) = 5$. The other entries of the altered energy function are computed analogously. Combining this with~(\ref{eq:cong_cond}), we complete the proof of Theorem~\ref{thm:G2identity}.

After principal specialisation, the grounded partitions corresponding to the paths in~\eqref{eq:paths} are
\begin{align*}
&0_h\,2_f\,7_e\,12_d\,17_c\,22_b\,24_h,\\
&0_h\,1_g\,6_f\,14_f\,19_e\,24_d\,30_b\,32_h,\\
&0_h\,1_g\,7_e\,15_e\,21_c\,26_b\,31_a.
\end{align*}

Using Theorem~\ref{thm:KMN2coloured}, equation~\eqref{eq:geq_grounded_char}, the same specialisation yields Theorem~\ref{thm:G2identity_geq}.

\section{Recursions}\label{section:recursions}

In this section, we express the sum side of the identity in Theorem~\ref{thm:G2identity_geq} as a system of recursions. We will use this to show that Theorem~\ref{thm:G2identity_geq} cannot be refined by keeping track of some of the colours in the generating functions as variables, as was done for Primc's identity by Dousse--Lovejoy~\cite{DL18} and by Dousse--Konan~\cite{DK22a}. To derive these recursions, we use Dousse's method of weighted words~\cite{D17}, which was motivated by Andrews' proof of Schur's identity~\cite{And67}.

First, we impose a total order on the set of colours $\mathcal{C}=\{a,b,c,d,e,f,g,h\}$, corresponding to the specialisations given in Figure~\ref{fig:specialisations}: $g > f > e > h > d > c > b > a$. This ordering, together with the natural total order on $\Z$, induces a total order on the set $\mathbb{Z}\times \mathcal{C}$ of coloured integers: $(k_1)_{c_1}>(k_2)_{c_2}$ if $k_1>k_2$ or if $k_1=k_2$ and $c_1>c_2$. For each colour $c \in \mathcal{C}$, consider the following generating functions:
\begin{align*}
P_{k_c} = P_{k_c}(t) := \sum_{\pi \in \mathcal{P}_{k_c}} C(\pi)t^{|\pi|} \qquad \text{and} \qquad R_{k_c} = R_{k_c}(t) := \sum_{\pi \in \mathcal{R}_{k_c}}C(\pi)t^{|\pi|},
\end{align*}
where \smash{$\mathcal{P}_{k_c}$} denotes the set of elements of \smash{$\mathcal{P}_{M}^{h,\geq}$} whose largest part, with respect to the total order~$>$ above, is $k_c$; and \smash{$\mathcal{R}_{k_c}$} denotes the set of elements of \smash{$\mathcal{P}_{M}^{h,\geq}$} whose largest part, with respect to the total order $>$ above, is at most $k_c$.

By definition, $R_{k_g} - R_{k_f} = P_{k_g}$, and similarly for other colours along the total order. First, compute $P_{k_g}$ directly by reading off column $g$ of the energy matrix shown in Figure~\ref{fig:G2energyfn}. The value in rows $a$, $b$, $c$, $d$ is $0$, so the parts $l_d$, $l_c$, $l_b$, $l_a$ can proceed $k_g$ for any $l\leq k$. In rows $e$, $f$, and~$h$, the value is $1$, so $k_e$, $k_f$, and $k_h$ cannot proceed $k_g$, but $l_e$, $l_f$ and $l_h$ can for any $l\leq k-1$. Finally, in row $g$, the value is $2$, so neither $k_g$ nor $(k-1)_g$ can proceed $k_g$, but $l_g$ can for $l\leq k-2$. We illustrate these observations below by painting the first part blue, the admissible parts that can proceed it are in green, and the forbidden parts are in red:
\begin{align*}
\textcolor{blue}{k_g} &> \textcolor{red}{k_f} > \textcolor{red}{k_e} > \textcolor{red}{k_h} > \textcolor{green!60!black}{k_d} > \textcolor{green!60!black}{k_c} > \textcolor{green!60!black}{k_b} > \textcolor{green!60!black}{k_a} \\
&> \textcolor{red}{(k - 1)_g} > \textcolor{green!60!black}{(k - 1)_f} > \textcolor{green!60!black}{(k - 1)_e} > \textcolor{green!60!black}{(k - 1)_h} \\
&> \textcolor{green!60!black}{(k - 1)_d} > \textcolor{green!60!black}{(k - 1)_c} > \textcolor{green!60!black}{(k - 1)_b} > \textcolor{green!60!black}{(k - 1)_a} > \cdots.
\end{align*}
Hence, we obtain the following recursion:
\[P_{k_g} = R_{k_g}-R_{k_f} = gt^k\bigl(P_{k_d} + P_{k_c} + P_{k_b} + P_{k_a} + R_{(k-1)_f}\bigr).\]
Analogously, we obtain the following six recursion:
\begin{align*}
 &P_{k_f}= R_{k_f}-R_{k_e} = ft^k\bigl(P_{k_b} + P_{k_a} + R_{(k-1)_e}\bigr), \\
 &P_{k_e}= R_{k_e}-R_{k_h} = et^k\bigl(P_{k_a} + R_{(k-1)_h}\bigr), \\
 &P_{k_d}= R_{k_d}-R_{k_c} = dt^k\bigl(P_{k_a} + R_{(k-1)_h}\bigr), \\
 &P_{k_c}= R_{k_c}-R_{k_b} = ct^k\bigl(P_{(k-1)_h} + R_{(k-1)_b}\bigr), \\
 &P_{k_b}= R_{k_b}-R_{k_a} = bt^k\bigl(P_{(k-1)_h} + R_{(k-1)_a}\bigr), \\
 &P_{k_a}= R_{k_a}-R_{(k-1)_g} = at^k\bigl(P_{(k-1)_h} + R_{(k-2)_g}\bigr).
\end{align*}
The only difference in the computation of $P_{k_h}$ is that row $h$ of column $h$ in the energy matrix has entry $0$, allowing the part $k_h$ to repeat indefinitely (which will be indicated by a green box around $k_h$ below). Consequently, we need to use a geometric series to express the generating function $P_{k_h}$ in terms of a generating function where the largest possible part is a strictly smaller coloured integer than $k_h$:
\begin{align*}
\textcolor{red}{k_g} &> \textcolor{red}{k_f} > \textcolor{red}{k_e} > \textcolor{green!60!black}{\boxed{\textcolor{blue}{k_h}}} > \textcolor{red}{k_d} > \textcolor{red}{k_c} > \textcolor{red}{k_b} > \textcolor{red}{k_a} \\
&> \textcolor{red}{(k - 1)_g} > \textcolor{red}{(k - 1)_f} > \textcolor{red}{(k - 1)_e} > \textcolor{green!60!black}{(k - 1)_h} \\
&> \textcolor{green!60!black}{(k - 1)_d} > \textcolor{green!60!black}{(k - 1)_c} > \textcolor{green!60!black}{(k - 1)_b} > \textcolor{green!60!black}{(k - 1)_a} > \cdots.
\end{align*}
Hence, the generating function $P_{k_h}$ satisfies
\begin{align*}
P_{k_h} &= R_{k_h}-R_{k_d} = ht^k\bigl(P_{k_h} + R_{(k-1)_g}\bigr),
\end{align*}
which can be written as
\[P_{k_h} = \frac{ht^k}{1-ht^k} R_{(k-1)_g}.\]

Ideally, we would like to refine Theorem~\ref{thm:G2identity_geq}, keeping some of the colours as variables. Since the coefficients are positive, this product expansion can only contain terms of the form
\[\bigl(-Xt^k;t^n\bigr)_\infty\qquad\text{or}\qquad\bigl(Xt^k;t^n\bigr)_\infty^{-1},\]
where $X$ is a monomial in the colour set $\mathcal{C}$ and $k$, $n$ are positive integers.

To achieve this, we use the recursions without any specialisations to compute the first few terms:
\begin{align*}
&1+ (a + b + c + d + e + f + g + h) t \\
&\qquad+ \bigl(a + b + c + d + a d + e + a e + f + (a + b) f + g + (a + b + c + d) g + h + h^2\bigr) t^2 \\
&\qquad+ \bigl(a + b + c + d + e + f + g + a d g + h + a h + h^3 + b (a + h) + c (a + b + h) \\
&\qquad\quad\ \, + d (a + b + c + d + h) + e (a + b + c + d + h) + f (a + b + c + d + e + h) \\
&\qquad\quad\ \, + g (a + b + c + d + e + f + h) + h (a + b + c + d + e + f + g + h) \bigr) t^3 + \cdots.
\end{align*}
From the coefficient of $t$, we can see that the infinite product must include either $(1+at)$ or~${1/(1-at)}$. Similarly, it must include either $(1+bt)$ or $1/(1-bt)$, and so on, for all colours~${a, b, c, \ldots, h}$. Therefore, the coefficient of $t^2$ must contain all \smash{$\binom{8}{2}$} pairwise products of these variables. However, this coefficient only has 17 monomial terms. Thus, it is impossible to keep any of the colours in the generating function without specialisation and still obtain an infinite product.

On the other hand, it might still be possible to keep some colours as variables in the dilated version, i.e., by mapping each colour to itself multiplied by its specialisation:
\begin{align*}
\begin{aligned}
 a &\mapsto aq^{-3}, \qquad b \mapsto bq^{-2}, \qquad c \mapsto cq^{-1}, \qquad d \mapsto d, \\
 h &\mapsto h, \qquad g \mapsto gq^{3}, \qquad f \mapsto aq^{2}, \qquad e \mapsto eq, \qquad t \mapsto q^{4}.
\end{aligned}
\end{align*}
We have attempted, via a case-by-case analysis, to keep some of the colours as independent variables in the dilated version, but without success. We believe that the resulting generating function does not admit a product form unless we specialise $a=b=\cdots=h=1$.

\subsection*{Acknowledgements}

I would like to express my gratitude to my PhD advisor, Jehanne Dousse, for suggesting this research question and for our weekly discussions on my progress. I greatly appreciate her generosity with her time and resources, as well as her valuable comments on earlier versions of this paper. I would also like to thank Matthew Russell for pointing out a mistake in an earlier version. I thank the anonymous referees for their valuable comments. The author is funded by the SNSF
Eccellenza grant of Jehanne Dousse, PCEFP2 202784.

\pdfbookmark[1]{References}{ref}
\LastPageEnding

\end{document}